\documentclass[a4paper,12pt]{article}
\usepackage{amsmath}
\usepackage{amssymb,amsfonts}
\usepackage[english]{babel}
\usepackage{gastex}
\usepackage{graphicx}
\usepackage{cite}
\usepackage[margin=3cm]{geometry}

\newtheorem{lemm}{Lemma}
\newtheorem{Rem}{Remark}
\newtheorem{prop}{Proposition}
\newtheorem{theor}{Theorem}

\begin{document}

\begin{flushleft}
{\large \bf ALMOST LIE NILPOTENT VARIETIES OF ASSOCIATIVE ALGEBRAS}

\vspace{1cm}

{\bf Olga~Finogenova}
\end{flushleft}


\noindent{\bf Abstract} 
We consider associative algebras over a field.
An algebra variety is said to be {\em Lie nilpotent}
if it satisfies a polynomial identity of the
kind $[x_1, x_2, \ldots, x_n] = 0$ where $[x_1,x_2] = x_1x_2 - x_2x_1$ and
$[x_1, x_2, \ldots, x_n]$ is defined inductively by
$$[x_1, x_2, \ldots, x_n]=[[x_1, x_2, \ldots, x_{n-1}],x_n].$$
By Zorn's Lemma every non-Lie nilpotent variety contains a minimal 
such variety, called 
{\em almost Lie nilpotent}, as a subvariety. 
A description of 
almost Lie nilpotent varieties for algebras over a field 
of characteristic 0 was made up by Yu.Mal'cev. 
We find a list of non-prime almost Lie nilpotent 
varieties of algebras over a field of positive characteristic.  


\section*{Introduction}

It is well-known that
every associative algebra may be considered as a Lie algebra under 
the Lie multiplication $[x,y] = xy-yx$. 
Structures of these algebras are closely connected, and properties
of one of them provide corresponding properties of other.
One of the most natural and intensively studied connections 
is an influence of Lie algebra nilpotency  on the associative structure.
The property means a fulfilment of the identity 
$[\cdots[[x_1,x_2],\ldots,x_n]=0$,
and the associative algebras (or varieties) satisfying the identity are called 
{\em Lie nilpotent}. Lie nilpotent algebras possess many wonderful
properties. For example, finitely generated such algebra is
finitely presented, residually finite, right and left Noetherian, 
representable by endomorphisms, etc. (cf.~\cite{Kublan}, ~\cite{Mark}).

Unfortunately, establishing Lie nilpotency  by means of explicit deriving of 
the required identity is, as a rule, a non-trivial and computationally 
laborious process. 
To facilitate it one can use a list of 
minimal w.r.t. inclusion elements
in the set of all non-Lie nilpotent varieties. 
We shall call such minimal varieties {\em almost~L.N.}.
By Zorn's Lemma 
every non-Lie nilpotent variety contains as a subvariety at least one 
almost~L.N. variety. Therefore, to prove Lie nilpotency of a variety
it is sufficient to verify that no one such varieties is contained
in our variety. This approach was successfully realized and corresponding 
descriptions were found out for many properties
such as commutativity~\cite{Mal-com},~\cite{Malc}, 
Engel condition~\cite{Engel}, 
locally residual finitness, locally weak Noetherian condition~\cite{Kublan}
and for many others. 

A description of almost~L.N. varieties of algebras over a field of 
characteristic zero is obtained in~\cite{Mal-net}. Note that 
the found there varieties have an 
elementary structure, due to which checking for being Lie nilpotent in this case 
is quite simple. The aim of the present paper is to make more clear a situation
with almost~ L.N. varieties for algebras over a field of positive characteristic. 

Hencerforth, all algebras are assumed associative.

We adopt the following notation.

Let $F$ be a field. We denote by $F\langle X\rangle$ the free
$F$-algebra generated by the countable set $X$. As usual, the elements of $F\langle X\rangle$
are called {\em polynomials}.
An ideal $I$ of $F\langle X\rangle$ is called {\em a $T$-ideal} if it is closed 
under endomorphisms. Let $A$ be an algebra, $\Sigma$ a set of polynomials, 
$\cal V$ a variety. We denote by ${\sf var} A$ the variety generated 
by $A$, and by 
${\sf var} \Sigma$ the variety defined by $\Sigma$.
An identity ideal of a variety ${\cal V}$ ( or an algebra $A$) is the set of 
all polynomials $f(x)$ such that $f(x)=0$ is an 
identity of ${\cal V}$( or $A$). We denote the 
identity ideal by $T({\cal V})$ (or $T(A)$, resp.);  and we write 
$T(\Sigma)$ in place of $T({\sf var} \Sigma)$.

We denote by ${\cal V}^{\star}$ the variety dual to ${\cal V}$; 
in the case of algebras $A^{\star}$ stands for the algebra anti-isomorphic 
to $A$.

By $\bar x$ we denote a tuple of variables
$x_1, x_2, \ldots \; $. It will always 
be clear from the context whether such is assumed ordered or not.

For a convenience we denote by  $W_n(\bar x)$ the Lie commutator of a
length $n$,  in other words, $W_2(x_1,x_2) = [x_1,x_2]$ and, defining inductively, 
$$W_n(x_1,\ldots, x_n) = [W_{n-1}(x_1,\ldots x_{n-1}), x_n].$$

\medskip

Let us introduce the notation needed for algebras. 

Denote by $C$ 
the algebra over a field generated by 
$c_1,c_2, \ldots $ with relations 
$$c_ic_j=c_jc_i, \;c_i^2=0,
\;i,j=1,2\ldots $$
If $p>0$ is the base field characteristic, then 
$C$ generates the variety  ${\sf var} \{[x,y]=0, x^p=0\}$.

Letting $U$ be an arbitrary
algebra, we put
$$A(U) \cong \left(
\begin{array}{cc}
U & U \\
0 & 0
\end{array}\right).$$

Recall that an algebra variety is said to be {\em prime} 
( or {\em verbally prime}) if for its T-ideal $T$
every inclusion  $I_1 \cdot I_2 \subseteq T$ holding for 
two T-ideals $I_1$ and $I_2$
implies either
$I_1 \subseteq T$ or $ I_2 \subseteq T$.

\medskip

\begin{theor}
Let F be an infinite field of a positive
characteristic, and let $\cal V$ be a non-prime F-algebra variety.
Then $\cal V$ is almost~ L.N. if and only if 
$\cal V$ is generated either by $A(C)$ or by $A(C)^{\star}$.

\end{theor}

\medskip

For a finite field $F$ we also introduce the following series of algebras:
$$B(F,G,\sigma) \cong
\left\{\left(
\begin{array}{cc}
b & c \\
0 & \sigma (b)
\end{array}\right)\right\}$$
where $b,c$ run through a finite extension $G$ of $F$,
and $\sigma$ is an automorphism of $G$ such that the
invariant field $G^{\sigma}$ is an unique maximal
subfield of $G$ containing $F$.

\begin{theor}

Let F be a finite field, and let
$\cal V$ be a non-prime F-algebra variety.
Then $\cal V$ is almost~ L.N. if and only if
$\cal V$ is generated by one of the following algebras
$A(F)$, $A(F)^{\star}$, $A(C)$, $A(C)^{\star}$, $B(F,G,\sigma)$.

\end{theor}

Although there is no information about prime almost~ L.N.
varieties in the descriptions, we think that sometimes 
the results may be quite effective. 
To demonstrate it we apply the descriptions
in a partial case (see {\bf Example}).

\section*{Almost~ L.N. varieties}

We prove the Theorems at the end of the section collecting some facts.

\begin{lemm}
\label{tech}
Let $n \ge 3$ be an integer and
$W_n(\bar x)=0$ be an identity of
a variety ${\cal M}$.
 Then
for every $k \le n-3$ 
${\cal M}$ satisfies
the identity
$$W_{n-k}(\bar y_1)\cdot W_{n-k}(\bar y_2) \cdots W_{n-k}(\bar y_{2^k}) = 0.$$
\end{lemm}

\noindent {\bf Proof.}
Let us substitute $x_{n-1}\mapsto zx_{n-1}$  into $W_n$.
We obtain  
$$0=[W_{n-2},zx_{n-1},x_n] $$
$$ = z[W_{n-2},x_{n-1},x_n] + [z,x_n][W_{n-2},x_{n-1}]$$
$$ + [W_{n-2},z][x_{n-1},x_n] +[W_{n-2},z,x_n] x_{n-1}.$$
The first and fourth summands belong to  $T(W_n)$.
Put $z=W_{n-2}(\bar u)$ and see that because of Jacobi identity 
the first commutator in the third summand becomes a consequence of  $W_n$.
Then we have modulo $T(W_n)$ 
$$[W_{n-2}(\bar u),x_n][W_{n-2}, x_{n-1}]=0,$$
that is $W_{n-1}(\bar y_1)\cdot W_{n-1}(\bar y_2)=0$. 
As far as $W_s u W_m \in T(W_s \cdot W_m )$ for arbitrary $s, m$ 
we can repeat the above argument  for 
each of $W_{n-1}$ to obtain a product of $4$ factors $W_{n-2}$ and so on.

\medskip

For the rest of this section we will denote by ${\cal V}$ an 
almost~ L.N. variety of algebras over a field (finite or infinite).
(Recall that  an almost~ L.N. variety is non-Lie nilpotent while
its proper subvarieties are all Lie nilpotent.)

\begin{lemm}
\label{simple}
Let $f(\bar x) \notin T({{\cal V}})$.
There exists an integer $n>0$ such that 
$$W_n(\bar y)  \in T(\{f\}) + T({\cal V}).$$
\end{lemm}

\noindent {\bf Proof.}
By the condition the ideal $ T(\{f\}) + T({\cal V})$ defines a proper
and, hence, Lie nilpotent  subvariety 
of ${\cal V}$.

\begin{lemm}
\label{useful}
If $T(\{f\})\cdot T(\{g\}) + T(\{g\})\cdot T(\{f\}) \subseteq T({\cal V})$ 
then either $f\in T({\cal V})$ or  $g(\bar x) \cdot g(\bar y) \in T({\cal V})$.
\end{lemm}

\noindent {\bf Proof.}
Arguing by contradiction,
suppose that $f\notin T({{\cal V}})$ and 
$g(\bar x) \cdot g(\bar y) \notin T({{\cal V}})$. By Lemma~\ref{simple}
we have an identity
$$W_n(\bar x) = h(\bar x)$$ 
for some $h \in T(\{f\})$. Furthermore,  
we have an identity
$$W_s(\bar y) = t(\bar y),$$ 
where $t(\bar y)$ is a sum of polynomials
of the following kind 
$$u g(a_1, a_2,  \ldots) g(b_1, b_2,  \ldots) v$$ 
and 
$u, v$ are monomials, $a_i, b_i$ are polynomials depending on $y_1,\ y_2\ldots$.
We can obtain the following equations modulo $T({\cal V})$:
$$W_{n+s-1}=W_s(W_n(\bar x), y_2, \ldots)=W_s(h(\bar x), y_2, \ldots)=t(h(\bar x), y_2, \ldots)$$
Now we show that $t(h(\bar x), y_2, \ldots) \in
T(\{f\})\cdot T(\{g\}) + T(\{g\})\cdot T(\{f\})$. 
The substitution $y_1 \mapsto h(\bar x)$ transforms every summand from $t$
in which $u$ or $v$ contains $y_1$ into a polynomial from 
$T(\{f\})\cdot T(\{g\}) + T(\{g\})\cdot T(\{f\})$. 
Consider another summands.
It is easy to see that 
$g(a_1, a_2,\ldots) = g(c_1, c_2, \ldots) + d(\bar y)$ and polynomials 
$c_1, c_2, \ldots$ do not depend on $y_1$ while every monomial of  $d(\bar y)$
contains $y_1$. 
Hence, the substitution $y_1 \mapsto h(\bar x)$ turns 
all summands of the kind 
$$u d(\bar y) g(b_1, \ldots)v$$ into a polynomial from
$T(\{f\})\cdot T(\{g\})$. 

In the other summands $ug(c_1, c_2, \ldots) g(b_1, b_2,  \ldots) v$
the variable $y_1$ has to occur in every monomial of $g(b_1, b_2,  \ldots)$.
Therefore, our substitution turns these summands into polynomials from 
$T(\{g\})\cdot T(\{f\})$.
Hence, $t(h(\bar x), y_2, \ldots) \in
T(\{f\})\cdot T(\{g\}) + T(\{g\})\cdot T(\{f\})$ and $W_{n+s-1} \in T({\cal V})$.
A contradiction.

\begin{lemm}
\label{comm}
If ${\cal V}$ is non-prime then it satisfies the identity
$[x,y][z,t]=0$.
\end{lemm}

\noindent {\bf Proof.}
As far as
${\cal V}$ is non-prime, there exist two polynomials $f_1, f_2 \notin T({\cal V})$ 
such that $T(\{f_1\})\cdot T(\{f_2\}) \subseteq T({\cal V})$. 
By Lemma~\ref{simple},
$W_n \in T(\{f_1\}) + T({{\cal V}})$ and $W_s \in T(\{f_2\}) + T({{\cal V}})$.
 By Lemma~\ref{tech},  we have for some $r$ 
$$\underbrace{W_3 \cdots W_3}_{r} \in 
(T(\{f_1\}) +T({\cal V})) (T(\{f_2\}) + T({{\cal V}})) \subseteq T({\cal V}).$$
Let $r$ be a minimal such number. 
Put $f = \underbrace{W_3\cdots W_3}_{r-1}$ and $g=W_3$.
Then, $f\notin T({\cal V})$ and 
$T(\{f\}) \cdot T(\{g\}) + T(\{g\}) \cdot T(\{f\}) \subseteq T({\cal V})$. 
By Lemma~\ref{useful},
$g(\bar x)\cdot g(\bar y)\in  T({\cal V})$. Therefore,
$W_3 W_3 = 0$ is an identity of ${\cal V}$.

Suppose,  $[x,y][z,t] \notin T({\cal V})$. By Lemma~\ref{simple},  ${\cal V}$ 
satisfies for some $m$ an identity                                     

\begin{equation}
\label{e111}
W_m(\bar x) = g(\bar x)
\end{equation}
where $g$ is a sum of polynomials of the kind
$u[x_i,x_j]v[x_l,x_k]w$ and $u,\,v,\, w$ are monomials depending on 
$x_1, \ldots, x_m$. It is easy to see that the substitutions
$x_i \mapsto [x_i,y_i]$ ($i=1,\ldots, m$) turn the right part of~(\ref{e111}) 
into
a consequence of the polynomial $W_3 \cdot W_3 \in T({\cal V})$. Therefore, the left 
part $W_m([x_1,y_1], \ldots ,[x_m,y_m])$ belongs to $T({\cal V})$ as well.
Let $m$ be a minimal number such that 
$$W_m([x_1,y_1], \ldots ,[x_m,y_m]) \in T({\cal V}).$$ 
We want to prove that
$m=2$. Suppose $m>2$.
By Lemma~\ref{simple}, for some $t$ the variety ${\cal V}$ 
satisfies  an identity 
\begin{equation}
\label{e2}
W_t(\bar x) = h(\bar x)
\end{equation}
where $h$ is a sum of polynomials of the kind
$$u W_{m-1}([a_1,b_1], \ldots ,[a_{m-1},b_{m-1}])v$$ and $u,\,v,\, a_i, \, b_i$ 
are monomials depending on 
$x_1, \ldots, x_m$. It is easy to see that  
$$[h(\bar x),[y,z]] \in T(\{W_m([x_1,y_1], \ldots ,[x_m,y_m])\}) + 
T(\{W_3(\bar x) W_3(\bar y)\}) \subseteq T({\cal V}).$$ 
Hence,
$[W_t(\bar x), [y,z]] = 0$ is an identity of ${\cal V}$. 
By the assumption, $m>2$, and
we have $t \ge 3$. 
Substitute $W_t(\bar y)$ for $x_1$ in $h(x_1, x_2, \ldots)$.
The summands of $h(\bar x)$ where $x_1$ occurs in $u$ or $v$ are transformed 
into polynomials from $T(\{W_3(\bar x) W_3(\bar y)\})$. Other summands are 
transformed into
sums of polynomials of
the kind $\tilde u [[c_1,c_2], c_3W_t(\bar y)c_4 ] \tilde v$ for some 
$\tilde u$, $c_1$, $c_2$, $c_3$, $c_4$, $\tilde v$. Clearly,
$$\tilde u [[c_1,c_2], c_3W_t(\bar y)c_4 ] \tilde v  \in 
 T(\{W_3(\bar x) W_3(\bar y)\}) + T(\{[W_t(\bar x), [y,z]]\}).$$
Hence, $h(W_t(\bar y), x_2, \ldots) \in T({\cal V})$.
Therefore,  replacing  $x_1$
by $W_t(\bar y)$ in~(\ref{e2}), we obtain an identity of
 ${\cal V}$ on the right and $$W_{2t-1} = W_{t}(W_{t}(\bar y), x_2, \ldots , x_t)$$
on the left. The contradiction shows that $m=2$, i.e. 
$[[x_1,y_1],[x_2,y_2]] = 0$ is an identity of ${\cal V}$. 

Substituting 
$y_2\mapsto y_2 z$,  we find the consequence 
$$[x_1,y_1,y_2][x_2,z] + [x_2,y_2][x_1,y_1,z]=0.$$ Substituting 
$y_2\mapsto y_2 t$ in the identity and using $W_3\cdot W_3 = 0$, we obtain
the consequence $[x_2,y_2][x_1,y_1,z,t] = 0$. 
By Lemma~\ref{useful}, we conclude that 
$[x_2,y_2][x_3,y_3] \in T({\cal V})$.

\medskip

The next property will help us to
simplify some identities. We say that a variety 
${\cal M}$ possesses {\bf \em Property Z} if
${\cal M}$ satisfies the following condition:

\smallskip

{\em for any polynomials $f$ and $g$,
if $f(h(\bar t), \bar x)\in T({\cal M})$ holds for all
$h \in T(\{g \})$, then either $g \in T({\cal M})$ or 
$f([y,z], \bar x)\in T({\cal M})$.}

\begin{lemm}
\label{z}
Let ${\cal V}$  satisfies an identity
$[x,y][z,t]=0$. Then ${\cal V}$  possesses 
Property Z.
\end{lemm}

\noindent {\bf Proof.}
Let $f$ and $g$ be polynomials from Property~Z condition.
Suppose that neither $g$ nor $f([y,z],\bar x)$ belongs to
$T({\cal V})$. By Lemma~\ref{simple}, ${\cal V}$ satisfies two identities
\begin{equation}
\label{e4}
W_r(\bar t) = h(\bar t),
\end{equation}
where $h(\bar t)$ is a consequence of $g$; and
\begin{equation}
\label{e3}
W_s(\bar u) = c(\bar u),
\end{equation}
where $c(\bar u)$ is a consequence of $f([y,z],\bar x)$.
Clearly, $c(\bar u)$ is a sum of summands of the two types:
\begin{itemize}
\item[A:] $vf(\sum_i [a_i,b_i] , \ldots)w$ and $u_1$ occurs in every 
monomial $a_i$,
\item[B:] $a[v,w]b$ and $u_1$ occurs either in $a$ or in $b$.
\end{itemize}
The substitution 
$u_1 \mapsto W_r(\bar t)$ in~(\ref{e3}) transforms
all summands of the type~B into consequences of $[x,y][z,t]$ and 
all summands of the type~A into polynomials from
$T(\{f(T(\{g\}),\ldots)\})$~(see~\ref{e4}).
Therefore, under this substitution we obtain modulo $T({\cal V})$
$$W_{r+s-1}=W_s(W_r(\bar t), u_2, \ldots ) = 0.$$  The contradiction
shows that either $g$ or      
$f([y,z], \bar x)$ belongs to $T({\cal V})$.
  
\medskip

In detail varieties satisfying  Property~Z are considered in~\cite{Engel}.
In the next item we give some necessary facts concerning such varieties.

Let $H$ be a relatively free countably generated  $F$-algebra of
$${\sf var}\{[x,y][z,t] = 0 \}.$$ Denote by $\Lambda$ the set of free generators
for $H$. For convenience, the elements of
$H$ will be called polymomials.

For an arbitrary monomial 
$f(x_1, x_2, \ldots )$ put $S_{x_i}(f)=\{m\}$ if $m$ is the number of
 occurences of the letter $x_i$
in $f$ (in other words, the degree of $f$ over $x_i$). 
If $f(\bar x, \bar t)= f_1(\bar x, \bar t) + \cdots + f_n(\bar x, \bar t)$ is
the sum of monomials $f_i$, we put
$$S_{\bar x}(f)=
\bigcup_{i=1}^n \bigcup_{u \in \bar x} S_u(f_i).$$

Similarly, for every polynomial $f$ in
$[H, H]$ we define 
$D(f)$, the set of {\it bilateral degrees}.
First, let $f(x, \bar t)$ be a commutator monomial of the form
$$a(x, \bar t)[t_i,t_j]b(x, \bar t)$$ with $S_x(a) = \{ k\}$ and
$S_x(b) = \{ l \}$. Then we put
$D_x(f) = \{(k,l)\}$.
If $f(x, \bar t)= a(x, \bar t)[x,t_i]b(x, \bar t)$ and $a, b$ are as above,
we put
$$D_x(f) = \{(k+1,s), (k,s+1)\}.$$ Finally, if 
$$f(\bar x, \bar t)= f_1(\bar x, \bar t) + \cdots + f_n(\bar x, \bar t)$$
is the sum of commutator monomials $f_i$,  we define
$D_{\bar x}(f)$, {\it the set of bilateral degrees of a polynomial  $f$ over 
variables $x_1, x_2, \ldots$}, as follows
$$D_{\bar x}(f) = \bigcup_{i=1}^n \bigcup_{u \in \bar x}D_{u}(f_i).$$

For instance, let 
$$f(x,y,z) = x[x,y]x^2 + y^5[y,z];$$ then
$S_{\{x,y\}}(f) =\{4, 1, 0, 6\}$
and $$D_{\{x,y\}}(f)=
\{(2,2), (1,3), (1,0), (0,1), (0,0), (6,0), (5,1)\}.$$

\begin{Rem} {\em Generally, the sets $S_{\bar x}(f)$ and $D_{\bar x}(f)$
are not  uniquely defined for  $f$ and depend on
its  particular representation --- as a sum of monomials or commutator
monomials. Below, by writing $S_{\bar x}(f) = M$ ($D_{\bar x}(f) = M$)
we mean that $f$ is given by a representation such that $S_{\bar x}(f)$ 
( $D_{\bar x}(f)$ resp.) coincides with the set $M$. 
Moreover, we assume that
$D_{\bar x}(f)$ and $S_{\bar x}(f)$ for a polynomial $f\in [H,H]$ are expressed
via the same representation, that is,
$S_{\bar x}(f)=\{ i_1+j_1, \ldots, i_m+j_m \}$ if
 $D_{\bar x}(f) =\{(i_1,j_1), \ldots, (i_m,j_m)\}$.}
\end{Rem}

Recall that polynomial $f(\bar x, \bar t)$ is said to be {\em essential} 
in variables from $\bar x$ if every $x_i$ occurs in every monomial of $f$.


\begin{prop}
\label{PropZ}
Let ${\cal M}$ be a $F$-algebra variety possessing  
Property Z and satisfying 
$[x,y][z,t]=0$. Suppose that
a polynomial $f(\bar x, \bar t) \in H$ essential in all 
$x_i$ has the form
$$ f(\bar x, \bar t) = w(\bar x, \bar t)
+ \sum_{i} g_i(\bar x, \bar t) + \sum_{i} h_i(\bar x, \bar t),$$
where $h_i(\bar x, \bar t)$ are monomials, $g_i(\bar x, \bar t)$ are commutator
monomials,\\ and  $w(\bar x, \bar t) \in [H,H]$. Also, assume that there exist 
finite sets
$$A \subseteq (\mathbb{N}\cup\{0\}) \times (\mathbb{N}\cup\{0\}) \mbox{ and }
B \subseteq \mathbb{N}\cup\{0\}$$
satisfying the following conditions:
\begin{itemize}
\item[1)] $A \cap D_{\bar x}(w(\bar x, \bar t)) = \emptyset\;$
and $\;B \cap S_{\bar x}(w(\bar x, \bar t)) = \emptyset;$
\item[2)] for every $i$, there exist
$j$ and $k$ such that\\
$D_{x_j}(g_i) \subseteq A\;$ and $\; S_{x_k}(h_i) \subseteq  B.$
\end{itemize}
If  $f \in T({\cal M})$ and $w(\bar x, \bar t) \notin T({\cal M})$ then
${\cal M}$ satisfyes an identity
$$x^l[y,z]x^m = x^r[y,z]x^s,$$
where
$(l,m) \in D_{ \bar x}(w(\bar x, \bar t))$ and
$(r,s) \in A$ or $r+s \in B$.

\medskip
(The set $A$, like $B$, may be empty, in which case
all $g_i$ (resp., $h_i$) will be zero.)  
\end{prop}

The proof of Proposition~\ref{PropZ} can be found in~\cite{Engel}.

Let us recall that by ${\cal V}$ we denote an almost~ L.N. variety.

\begin{lemm}
\label{glav}
Let ${\cal V}$  be generated by nilpotent  algebras and $[x,y][z,t] \in T({\cal V})$.
Then ${\cal V}$ satysfies either 
the identity $x[y,z]=0$ or the identity $[y,z]x=0$.
\end{lemm}

\noindent {\bf Proof.}
Suppose that $x[y,z] \notin T({\cal V})$. Then, by Lemma~\ref{simple}, 
${\cal V}$ satisfies 
for some $n$ an identity 
$$W_n(\bar x) = h(\bar x),$$
where $h$ is a sum of polynomials of the kind
$u x_i[x_j,x_k]v$ and $u,\,v$ 
are monomials (maybe, empty). Substituting $[y,z]$ for $x_1$ 
we transform
all summands where  $x_1$ does not occur
in commutators into consequences of $[x,y][z,t]$.  
Thus, under this substitution
we get an identity 
$$[y,z]x_2x_3\cdots x_n + g(y,z,\bar x) = 0,$$
where $g=\sum\limits_{i>1,j} a_{ij} x_i[y,z]b_{ij}$. By
Lemma~\ref{z},
${\cal V}$ possesses Property Z and satysfies all conditions of 
Proposition~\ref{PropZ}. Put 
$$w = [y,z]x_2x_3\cdots x_n,$$
$$A=\{(n,k)| n\ge 1\} \cap D_{\{x_2, x_3, \ldots, x_n\}}(g), \;\;\;B=\emptyset. $$ 
Clearly,
$D_{\{x_2, \ldots, x_n\}} (w) =\{(0,1)\}$ 
and $D_{\{x_2, \ldots, x_n\}} (w) \cap A =\emptyset $. 
It is easy to see that our $w$, $A$ and  $B$ satisfy  the conditions 
of Proposition~\ref{PropZ}.
Thus, ${\cal V}$ satysfies either the identity
$w =0$ or an identity 
$[y,z]x = x^r[y,z]x^s$ where $r,s$ are integer and $r\ge 1$.

Suppose, at first,  that 
$w=0$. We can assume that $n$ is a minimal  such number, that is 
$[y,z]x_2x_3\cdots x_{n-1}  \notin T({{\cal V}})$. 
Put $f(t,x) = tx$ and
$g=[y,z]x_2x_3\cdots x_{n-1}$. Clearly, $f(h,x) \in T({\cal V})$ for every
$h\in T(\{g\})$. Hence, we have $f([y,z],x) \in T({{\cal V}})$
because ${\cal V}$ satysfies the property~Z
and $g \notin T({{\cal V}})$. Thus, 
$[y,z]x \in T({{\cal V}})$. In this case
Lemma~\ref{glav} is proved.

Now suppose that $[y,z]x = x^r[y,z]x^s$ is an identity of ${\cal V}$.
We can assume $s=0$; otherwise, we obtain the identity $[y,z]x=0$
 because ${\cal V}$ is generated by nilpotent algebras.
Then, we see that $r \ne 1$, because $[[y,z],x] \notin T({\cal V})$.
Thus, an identity $[y,z]x = x^r[y,z]$ with $r>1$ holds in ${\cal V}$.
 
Furthermore,  assuming, in addition, that $[y,z]x \notin T({\cal V})$
we find by the dual argument an identity $x[y,z] = [y,z]x^k$ for an integer 
$k>1$. 
In every nilpotent 
algebra two identities $[y,z]x = x^r[y,z]$ and $x[y,z] = [y,z]x^k$ imply
$x[y,z] =0$ and $[y,z]x=0$ which contradicts the assumption 
because ${\cal V}$ is generated by nilpotent algebras.

\medskip

Below, $p>0$ is the characteristic of the base field.

\begin{lemm}
\label{yp}
If $[y,z]x \in T({\cal V})$ 
and ${\cal V}$  is generated by nilpotent  algebras then $y^px \in T({\cal V})$. 
\end{lemm}

\noindent {\bf Proof.}
It is easy to see that modulo $T(\{[y,z]x\})$  
the polynomial $W_n(\bar x)$ is equal to $x_3\cdots x_n[x_1,x_2]$
and $(u+v)^px$ is equal to $u^px+v^px$. 
Suppose that  $y^px \notin T({\cal V})$. 
By Lemma~\ref{simple} and  the above observation,  
${\cal V}$ satisfies an identity
$$x_3\cdots x_n [x_1,x_2]= h(\bar x)$$
where $h(\bar x)$ is a sum 
of monomials of the kind
$ab^p c$. Clearly, the length of every such monomial  
is greater than $n$.
Substituting $[x_1,y]$ for $x_1$ in the identity, we obtain
modulo $[y,z]x$
$$x_2x_3\cdots x_n[x_1,y] = \sum\limits_i u_i[y,x_1],$$
where the length of every monomial $u_i$ is greater than $n-1$.
Since ${\cal V}$ is generated by nilpotent algebras,
it satisfies the identity $ x_2x_3\cdots x_n[x_1,y] = 0$ and, hence,
$W_{n+1}(y,x_1, \ldots, x_n) = 0$. A contradiction. 
  
\medskip

Let us denote by ${\cal A}_p$ the variety given by 
$[y,z]x=0$ and $y^px=0$. Obviously, 
${\cal A}_p^{\star}$ is given by 
$x[y,z]=0$ and $xy^p=0$. 
Consider these varieties. First of all, they are not Lie nilpotent because of 
a following easily verified lemma.

\begin{lemm}
\label{a}
The algebras $A(C)$ and $A(C)^{\star}$ are not Lie nilpotent. 
The algebra $A(C)$  belongs to ${\cal A}_p$, and $A(C)^{\star}$  belongs to 
${\cal A}_p^{\star}$.  
\end{lemm}

Now, to show that ${\cal A}_p$ and ${\cal A}_p^{\star}$ are almost~L.N  
it suffices to prove 
the following lemma and its dual analogue.

\begin{lemm}
\label{proper}
Every proper subvariety of ${\cal A}_p$ is Lie nilpotent. 
\end{lemm}

\noindent {\bf Proof.}
Let ${\cal M}$ be a proper subvariety of ${\cal A}_p$. 
Then, there exists  
a polynomial $f$ from $T({\cal M}) \setminus T({\cal A}_p^{\star})$.
We want to find a consequence of $f$ of the form
\begin{equation}
\label{e6}
u_1^{s_1} \cdots u_m^{s_m}[y,z], \;\;\; s_1 < p, \ldots , s_m < p.
\end{equation}
Obviously, its full linearization and 
the polynomial  $[y,z]x$ have as a consequence the desired polynomial
$W_n(\bar x)$ for $n = (p-1)m+2$.

Since $[u,v]x_1+[x_1,u]v + [v,x_1]u = x_1[u,v]+v[x_1,u] + u[v,x_1]$ 
in every associative algebra,
we have modulo $T(\{[y,z]x\})$ 
$$x_1[u,v] = -v[x_1,u] - u[v,x_1].$$
Moreover,  we have modulo $ T({\cal A}_p^{\star})$ 
$$x_1^{p-1}[x_1,u] = ux_1^p.$$ 
Therefore,
$f$ can be written
modulo $ T({\cal A}_p)$ in a following   
form 

\begin{eqnarray*}
f(\bar x)=&
\sum\limits_{\bar k: \, k_1 < p, \ldots, k_n < p} &
\alpha_{\bar k} x_1^{k_1}x_2^{k_2}\cdots x_n^{k_n} + \\
&\sum\limits_{i,\, \bar s: \, s_1 < p-1, s_2<p , \ldots, s_n < p}& 
\beta_{i,\bar s} 
x_1^{s_1}x_2^{s_2}\cdots x_n^{s_n}[x_1,x_i] +\\
&\sum\limits_{\bar t: \, t_2 < p, \ldots, t_n < p}& \gamma_{\bar t } 
x_2^{t_2}\cdots x_n^{t_n} x_1^p 
\end{eqnarray*}

If for some $\bar k$ we have $\alpha_{\bar k} \ne 0$, we multiply $f$  by $[y,z]$ 
from the left
and obtain modulo $T({\cal A}_p)$
$$\sum\limits_{\bar k: \, k_1 < p, \ldots, k_n < p} 
\alpha_{\bar k}x_1^{k_1}x_2^{k_2}\cdots x_n^{k_n}[y,z] \in T({\cal M}).$$
Since degrees of variables are smaller than $p$, every polyhomogeneous
summand of the last polynomial belongs to $T({\cal M})$. It remains to see
that all of them are of the form~(\ref{e6}).

If $\alpha_{\bar k} =0$ for all $\bar k$ and $\beta_{i,\bar s} \ne 0$
for some pair $\bar s$, $i$, we 
substitute $x_1 +[y,z]$ for $x_1$ in $f$ and select all summands 
with $[y,z]$ to obtain modulo $T({\cal A})$
$$\sum\limits_{\bar s: \, s_1 < p-1, s_2<p , \ldots, s_n < p} \beta_{i,\bar s} 
x_1^{s_1+1}x_2^{s_2}\cdots x_n^{s_n}[y,z] \in T({\cal M}).$$
Repeating the above argument we obtain the desired 
polynomial~(\ref{e6}).

Finally, if all $\alpha$ and $\beta$ equal zero, we 
substitute $x_1 +[y,z]$ for $x_1$ in $f$ and select all summands 
with $[y,z]$ to obtain modulo $ T({\cal A})$
$$\sum\limits_{\bar t, t_2 < p, \ldots, t_n < p} \gamma_{\bar t } 
x_1^{p-1} x_2^{t_2}\cdots x_n^{t_n}[y,z]\in T({\cal M}).$$
As above, every polyhomogeneous
summand of the polynomial belongs to $T({\cal M})$ and it has the form~(\ref{e6}). 

Hence, every proper subvariety of ${\cal A}_p$ is Lie nilpotent.
  
\medskip

Now we are ready to prove the first main result.

\medskip

{\bf Proof of Theorem 1.}

Necessity. 
 Let ${\cal V}$ be a non-prime almost~ L.N. variety. 
It is well known that every variety of algebras over an infinite
field is generated by its nilpotent algebras. 
Hence, by Lemmas~\ref{comm},~\ref{glav},~\ref{yp} and dual to~\ref{yp},  
${\cal V} $ is contained in  ${\cal A}_p$ or in ${\cal A}_p^{\star}$.
By Lemmas~\ref{a} and~\ref{proper} and by their dual variants
${\cal A}_p$ and ${\cal A}_p^{\star}$ are almost~L.N. Hence, 
${\cal V}$ coincides with one of them. It remains to see that 
an almost~L.N. variety is generated by any its non L.N. algebra.
By Lemma~\ref{a}   
we can give $A(C)$ and $A(C)^{\star}$ as such
algebras. 

Sufficiency follows from  Lemmas~\ref{a},~\ref{proper} and their
dual analogues. 

\medskip

To prove  Theorem~2 we need in the discription of almost non-Engel varieties. Recall 
that a variety is said to be Engel
if for some natural $n$ it satisfies an identity $W_{n+1}(x,y,\ldots,y)=0$.
A variety is called {\em almost~ Engel} if it is itself non-Engel
but its proper subvarieties are all Engel.

\begin{prop} 
\label{justEn}
(Theorem 2, ~\cite{Engel})
A variety of algebras over a finite field $F$ 
is almost~Engel if and only if it is generated by one of the
algebras $A(F)$, $A(F)^{\star}$, or $B(F,G,\sigma)$.
\end{prop}

\begin{Rem} {\em In fact, every algebra in the list 
of Proposition~\ref{justEn} generates 
a non-Engel almost~ commutative variety (it is easy to 
verify immediately or see~\cite{Malc}). Hence, 
it generates an almost~ L.N. variety.}
\end{Rem}

\medskip

{\bf Proof of Theorem 2.}

Necessity. Let ${\cal V}$ be almost~ L.N.
Assume, at first, that ${\cal V}$ is non-Engel. 
Then, by Zorn's Lemma, it contains some almost Engel variety as a subvariety. 
Hence, by Proposition~\ref{justEn} and Remark~2 ${\cal V}$ coincides with one of
the varieties ${\sf var} A(F)$, ${\sf var} A(F)^{\star}$, or ${\sf var} B(F,G,\sigma)$.

Assume now, that  ${\cal V}$ is non-prime and Engel. 
It remains to show that in this case ${\cal V}$ is generated by nilpotent algebras.
Then, we shall be able to
repeat word for word the argument of Theorem~1 to
verify that ${\cal V}$ is generated 
by $A(C)$ or $A(C)^{\star}$.
So, being Engel ${\cal V}$ satisfies an identity 
$W_{n+1}(x,y,\ldots,y)=0$
Without loss of generality, we can assume that $n=p^t$ where $p$ is the
characteristic of the base field.
It is easy to prove that for every $n$
$$W_{n+1}(x,y,y,\ldots,y) = \sum\limits_{k=0}^n C_n^k (-1)^k y^k x y^{n-k}.$$
Therefore, for $n=p^t$ we have
$$W_{n+1}(x,y,y,\ldots,y) = [x,y^n].$$

Thus, ${\cal V}$ satisfies the identity $[x,y^{p^t}]=0$. 
Hence, the center of 
every algebra $A\in {\cal V}$  contains all subalgebras of $A$ which 
are finite fields. Moreover, 
${\cal V}$ is locally residual finite (see~\cite{Kublan}). Hence,
it is generated by its finite-dimensional algebras. Consider
an arbitrary such algebra $A$. 
The variety ${\cal V}$ does not contain full matrix algebras because each
of them is not Engel.
Hence,  $A=B \dot + J(A)$ and
$B$ is a finite direct sum of finite fields, $J(A)$ is a nilpotent 
radical. By above remark, $B$ is contained in the center of $A$. Therefore,
$W_k(A,A, \ldots , A) \subseteq W_k(J(A),J(A), \ldots , J(A))$ for every $k$.

Suppose now that  all nilpotent algebras of ${\cal V}$ generate a proper subvariety.
Since ${\cal V}$ is an almost~L.N. variety, this proper 
subvariety satisfies an identity
$W_k(\bar x)=0$. This identity holds in all algebras $J(A)$. 
Therefore, by the last inclusion  it also holds in all 
finite-dimensional algebras $A$.
Thus, $W_k(\bar x)=0$ is an identity of ${\cal V}$. The contradiction
shows that ${\cal V}$ is generated by nilpotent algebras.

The sufficiency follows from  Remark~2, Lemmas~\ref{a},~\ref{proper} and their
dual analogues.

\bigskip

We found descriptions of almost~ L.N. varieties  modulo prime varieties.
Nevertheless, we hope to show below that the descriptions are quite useful.

\section*{Example}

Let $p$ be the characteristic of a field $F$. Denote by 
${\cal M}$ the variety of $F$-algebras given by two identities
\begin{equation}
\label{comp}
[x^p,y]=0, 
\end{equation}
and 
\begin{equation}
\label{nilp}
W_{n_1}(\bar x_1)\cdot W_{n_2}(\bar x_2) \cdots W_{n_s}(\bar x_s) = 0
\end{equation}
for some natural $n_1, \ldots, n_s$.
We state that ${\cal M}$ is Lie nilpotent.

Indeed, assume, on the contrary, that ${\cal M}$ is not Lie nilpotent.
By Zorn's Lemma ${\cal M}$ contains an almost~L.N. variety ${\cal N}$ as a subvariety.
The subvariety is not verbally prime because of~(\ref{nilp}).
Hence, ${\cal N}$  is generated by  one of the algebras from Theorem lists.
Algebras $A(F)$, $A(F)^{\star}$, $B(F,G,\sigma)$ from Theorem~2 are
non-Engel (see Proposition~\ref{justEn}). Therefore, they do not satisfy 
Identity~(\ref{comp}), that is $p$-Engel condition, 
and can not belong to ${\cal M}$. 

For $A(C)$
put 
$x= \left(
\begin{array}{cc}
c_1 +\cdots + c_{p-1} & c_p \\
0 & 0
\end{array}\right)$,
$y = \left(
\begin{array}{cc}
c_{p+1} & 0 \\
0 & 0
\end{array}\right)$.
It is easy to prove that 
$x^p= \left(
\begin{array}{cc}
0 & (p-1)!\cdot c_1\cdots c_{p-1} c_p \\
0 & 0
\end{array}\right)$,
and 
$$[x^p,y]= x^py-yx^p = \left(
\begin{array}{cc}
0 & 0 \\
0 & 0
\end{array}\right) -
\left(
\begin{array}{cc}
0 & (p-1)!\cdot c_1\cdots c_{p+1}\\
0 & 0
\end{array}\right) \ne 0.$$
Thus, $A(C)$ and $A(C)^{\star}$ (by dual reason) 
also do not belong to ${\cal M}$. None of the algebras can generate the variety
${\cal N}$. The contradiction shows that ${\cal M}$ is Lie nilpotent.

\end{document}